\documentclass[reqno,11pt]{amsart}
\usepackage{graphicx}
\usepackage{amssymb}
\usepackage{amsmath}
\usepackage{amsfonts}
\usepackage{amsthm}
\usepackage{color}
\usepackage{hyperref}
\usepackage[mathscr]{eucal}
\usepackage{enumerate}
\usepackage{times}
\usepackage{authblk}
\usepackage[foot]{amsaddr}

\numberwithin{equation}{section}
\allowdisplaybreaks[1]

\newtheorem{Def}{Definition}[section]
\newtheorem{Thm}[Def]{Theorem}

\newtheorem{Lemma}[Def]{Lemma}

\newcommand{\beq}{\begin{equation}}
\newcommand{\eeq}{\end{equation}}
\newcommand{\Proof}{\begin{proof}}
\newcommand{\QED}{\end{proof} \noindent}

\newcommand{\mm}{\hspace{-.08cm}\cdot \hspace{-.08cm}}

\newcommand{\M}{\mathcal{M}}

\newcommand{\R}{\mathbb{R}}

\newcommand{\Gammati}{\tilde{\Gamma}}
\newcommand{\Riem}{{\rm Riem}}

\title[Geodesic Equation with Curvature Bounds]{On Weak Solutions to the Geodesic Equation\\ in the Presence of Curvature Bounds}

\author[M.\ Reintjes]{Moritz Reintjes$^*$}
\address[*]{Department of Mathematics\\ City University of Hong Kong \\ SAR Hong Kong}
\email{moritzreintjes@gmail.com}

\author[B.\ Temple]{Blake Temple$^{**}$ \\ \\ June 8, 2023}
\address[**]{Department of Mathematics\\ University of California\\ Davis, CA 95616\\ USA}
 \email{temple@math.ucdavis.edu}

\begin{document}

\begin{abstract}
We show that taking account of bounded curvature reduces the threshold regularity of connection coefficients required for existence and uniqueness of solutions to the geodesic equation, to $L^p_\text{loc}$, one derivative below the regularity $W^{1,p}_\text{loc}$ required if one does not take account of curvature, ($L_{loc}^p$ for existence, $W_{loc}^{1,p}$ for uniqueness). Our argument is based on authors' theory of the RT-equations for regularizing connections to optimal regularity by coordinate transformation. The incoming regularity is too low to formulate a weak version of the geodesic equation based on the standard method of multiplying by smooth test functions and integrating by parts, so alternatively, we define weak solutions by coordinate transformation and we give an explicit procedure for mollifying the original connection such that the correct weak solution is indeed a limit of smooth solutions of the mollified equations in the original coordinates.   This is an example where limits under suitable mollification are more fundamental than a weak formulation, indicative of more complicated PDE's in which the standard weak formulation of the equations does not adequately rule out unphysical solutions. Our results apply to general second order ODE's in which the lack of regularity can be isolated in the connection coefficients. The results apply to General Relativity.                                             
\end{abstract}

\maketitle 


\section{Introduction}  \label{Sec_intro}

We introduce a solvability condition sufficient to imply existence and uniqueness of solutions $x=\gamma(t)$ to the initial value problem for the geodesic equation\footnote{We use standard tensor notation;  indices $\mu,\nu,\rho,...$ run from $1$ to $n$, repeated up-down indices are summed over, etc., (see for example \cite{HawkingEllis}).}    
\beq \label{geodesic_eqn}
\begin{cases}
\ddot\gamma^\mu + \Gamma^\mu_{\rho\nu}(\gamma) \dot\gamma^\rho \dot\gamma^\nu =0,  \cr
\gamma^\mu(t_0)=x^\mu_0,  \cr 
\dot\gamma^\mu(t_0) = v^\mu_0,
\end{cases}
\eeq
when connection components are only in $L^p$, one derivative less regular than the standard ODE theory requires---a regularity too low to restrict connections to curves and make sense of weak solutions of \eqref{geodesic_eqn} in a standard way.  Our improvement is obtained by writing the equations in coordinates where the coefficients are more regular. For this, based on authors' prior work \cite{ReintjesTemple_ell2, ReintjesTemple_ell4}, it suffices to take account of the regularity of components of the Riemann curvature, $\Riem(\Gamma) \in L^p$ for existence and $\Riem(\Gamma) \in W^{1,p}$ for uniqueness of solutions to \eqref{geodesic_eqn}, and no other {\it apriori} information about the geometry associated with $\Gamma$ need be assumed.\footnote{Such as positive definiteness, c.f. \cite{KazdanDeTurck, Uhlenbeck}, or the need for $\Gamma$ to be a metric connection.}  To start, assume that the connection components $\Gamma \equiv \Gamma_x \equiv   \Gamma^\mu_{\rho\nu}(x) \in L^p(\Omega)$ are arbitrary given real valued functions of $x = (x^1,...,x^n)\in\Omega\subset \R^n,$  $\Omega$ open, $n\geq 2$. At this level of generality, the component functions $\Gamma^\mu_{\nu\rho}(x)$, together with an atlas of coordinate transformations and the transformation law for connections, are sufficient to define a unique affine connection $\Gamma$ on the tangent bundle of an $n$-dimensional manifold $\M$--and even though the associated geometry could be non-metric and highly degenerate, this alone is sufficient to define the Riemann curvature tensor $\Riem(\Gamma)$ associated with $\Gamma$, to which our theory here applies.  Thus, from the point of view of ODE theory, we can interpret the geometry as a device for formulating a solvability condition for general systems of nonlinear equations of form (\ref{geodesic_eqn}), whether or not the underlying geometry is of interest in its own right.

Note that equation \eqref{geodesic_eqn} does not admit a weak formulation based on multiplying by smooth test functions and integrating by parts because the connection $\Gamma\in L^p$ is of too low a regularity to restrict to curves $\gamma(t)$.  We here introduce an alternative formulation of the equations based on coordinate transformation, which is equivalent for smooth $\Gamma$, and hence suffices as a weak formulation when $\Gamma$ is of low regularity. The idea is that for smooth $\Gamma$, whenever we transform a solution $\gamma_x(t)$ of \eqref{geodesic_eqn} as a curve under a smooth coordinate transformation $x\to y$, the theory of the covariant derivative implies that to recover the equivalent equation for solution $\gamma_y(t)$ in the new $y$-coordinates, it is sufficient for the connection coefficients $\Gamma^\mu_{\rho\nu}(x)$ to transform by the connection transformation law.  By this principle of equivalent equations, the transformation law for connections is invoked by the requirement that the transformed equations be equivalent for smooth enough $\Gamma$, independent of any geometry in the background. It follows that if we can find a coordinate transformation that sufficiently regularizes $\Gamma$, then the transformed equations can naturally be taken as providing the correct weak formulation when the untransformed $\Gamma$ has too low a regularity to make sense of the equations in the untransformed coordinates. That is, if a low regularity $\Gamma\in L^p$ admits a regularization under coordinate transformation $x\to y$ sufficient to define classical solutions $\gamma_y(t)$ in the transformed coordinates, (our purpose is to establish this here using the theory of the RT-equations \cite{ReintjesTemple_ell4, ReintjesTemple_ell6}), then transforming the solution $\gamma_y(t)$ back to the original $x$-coordinates as a curve, $\gamma_x=y^{-1}\circ\gamma_y$, provides the correct notion of weak solution in the original coordinates.    Thus, analogous to defining distributions, we use two formulations equivalent for smooth solutions to define weak solutions in a low regularity setting in which one formulation, but not the other, has sufficient regularity to define classical solutions.   This principle works to define weak solutions for general second order systems of ODE's which are quadratic in first derivatives in the terms with coefficients of low regularity, i.e., equations of the form
\beq \label{gen_geodesic_eqn}
\ddot\gamma^\mu + \Gamma^\mu_{\rho\nu}(\gamma) \dot\gamma^\rho \dot\gamma^\nu =K^{\mu}(t,\gamma,\dot{\gamma}),
\eeq
where $K^\mu$ encodes the sufficiently regular terms in the equation.\footnote{Equation \eqref{gen_geodesic_eqn} is the relativistic version of Newton's equation of motion with $K$ the external Minkowski-force, \cite{Weinberg}.} In this case, for smooth $\Gamma$ and $K$, a smooth transformation of coordinates will transform $\Gamma$ by the connection transformation law, and $K^{\mu}$ must transform as a vector, $K^i=K^{\mu}\frac{\partial y^i}{\partial x^\mu}$, in order to get an equivalent equation of form \eqref{gen_geodesic_eqn} in $y$-coordinates. To see this note that the first two terms in \eqref{gen_geodesic_eqn} comprise the covariant derivative of $\nabla_{\dot{\gamma}} \dot{\gamma}$ which transforms as a vector, so transforming $K$ as a vector suffices to make \eqref{gen_geodesic_eqn} a covariant geometric equation.   Thus when $\Gamma\in L^p$, if the map $x\to y$ and $K^{\mu}$ have sufficient regularity for the equations (Lipschitz continuity suffices) to admit strong solutions in $y$-coordinates, then this will provide the correct weak formulation in $x$-coordinates. To keep things simple here, we restrict to the case $K^{\mu} = 0$.

It is well known that Peano's Theorem establishes H\"older continuity of the components $\Gamma\equiv \Gamma^\mu_{\nu\rho}(x) \in C^{0,\alpha}$, $0\leq \alpha < 1$, as the threshold regularities for existence of classical solutions of \eqref{geodesic_eqn} in $x$-coordinates.\footnote{Regularity of connections and tensors always refers to regularity of their components in a given coordinate system.} Uniqueness requires more regularity. Namely, the Picard-Lindel\"off Theorem requires Lipschitz continuity of the connection components, $\Gamma \in C^{0,1}$. By Morrey's inequality, functions in $W^{1,p}$, (the Sobolov space of functions with weak derivatives integrable in $L^p$), are H\"older continuous when $p>n$, $W^{1,p} \subset C^{0,\alpha}$ with $\alpha = 1- n/p$, and by Rademacher's Theorem $C^{0,1} \simeq W^{1,\infty}$, c.f. \cite{Evans}. So $\Gamma \in W^{1,p}$ is a regularity sufficient for uniqueness of solutions to \eqref{geodesic_eqn} and $\Gamma \in W^{1,\infty}$ is the threshold for uniqueness.   Below these threshold regularities, there is no guarantee the Picard iteration for solving (\ref{geodesic_eqn}) will converge to a solution, and prior approaches \cite{SaemannSteinbauer, Steinbauer} relied on the Filippov theory for establishing existence of (non-unique) solutions in the sense of differential inclusions \cite{Filippov}. In this paper, for regularities below these thresholds, we establish convergence of the Picard iteration in modified coordinates, under natural conditions on the regularity of $d\Gamma$, or equivalent, ${\rm Riem}(\Gamma)$.  

As authors pointed out in \cite{ReintjesTemple_ell3}, a direct consequence of Riemann's idea to construct a tensorial measure of curvature is that the regularity of connections can range by coordinate transformation from one derivative above to one derivative below the regularity of the Riemann curvature, while keeping the regularity of the Riemann curvature fixed. This is based on the difference between the transformation law for components of the Riemann curvature tensor $R^\tau_{\mu\nu \rho}(x) \equiv {\rm Riem}(\Gamma_x)$,
\beq \label{Riemann_transfo}
R^\tau_{\mu\nu \rho}(x) 
= \tfrac{\partial x^\tau }{\partial y^\delta} \; \tfrac{\partial y^\alpha}{\partial x^\mu } \tfrac{\partial y^\beta}{\partial x^\nu } \tfrac{\partial y^\gamma}{\partial x^\rho } \; R^\delta_{\alpha \beta \gamma}(y),
\eeq
and the transformation law for connections
\beq \label{connection_transfo}
(\Gamma_x)^\mu_{\rho\nu} = \tfrac{\partial x^\mu}{\partial y^\alpha } \Big( \tfrac{\partial y^\beta}{\partial x^\rho} \tfrac{\partial y^\gamma}{\partial x^\nu }  \, (\Gamma_y)^\alpha_{\beta \gamma}   +  \tfrac{\partial }{\partial x^\rho}  \tfrac{\partial y^\alpha}{\partial x^\nu } \Big),
\eeq
under a coordinate transformation $x^\mu\to y^\alpha$ with Jacobian $J \equiv \frac{\partial y^\alpha}{\partial x^\mu}$ and inverse $J^{-1} \equiv \frac{\partial x^\mu}{\partial y^\alpha}$, where $\Gamma_x$ and $\Gamma_y$ denote connection components represented in $x$- and $y$-coordinates, respectively, (e.g., $\Gamma_x\equiv (\Gamma_x)^\mu_{\rho\nu} \equiv (\Gamma_x)^\mu_{\rho\nu}(x)$, etc.). It follows from \eqref{Riemann_transfo} and \eqref{connection_transfo} that the regularity of connections and tensors, as measured in coordinates, is invariant under smooth coordinate transformations, but it is not invariant under low regularity transformations. In particular, the regularity of the components of the Riemann curvature tensor is invariant under the atlas of coordinate transformations with Jacobian of (at least) the regularity of the curvature, because the curvature components transform as a tensor by contraction with undifferentiated Jacobians. But, within the same atlas, the regularity of connection components is {\it not} invariant, due to the fact that \eqref{connection_transfo} involves derivatives of the Jacobian. By this, in a given coordinate system, the regularity of connection components could be up to one derivative {\it below} the regularity of the components of its Riemann curvature.   On the other hand,  the authors' theory of the RT-equation asserts that one can always transform to coordinates in which connection components are one derivative more regular than curvature components, (optimal regularity, c.f. \cite{ReintjesTemple_ell2, ReintjesTemple_ell4, ReintjesTemple_ell5}).  Thus, in a given coordinate system, the regularity of connection components could range from one derivative above, to one derivative below, the regularity of the components of its Riemann curvature.  It follows that, {\it apriori}, there is no {\it geometric} reason why the connection regularity in \eqref{geodesic_eqn} could not lie anywhere within this range.  We show in this paper that by taking account of the regularity of the curvature, the author's theory of optimal regularity based on the {\it Regularity Transformation (RT)-equations} \cite{ReintjesTemple_ell2, ReintjesTemple_ell4, ReintjesTemple_ell5}  can be applied to lower the threshold required for existence and uniqueness of \eqref{geodesic_eqn} by one derivative, over and above that given by the standard theory of Peano and Picard-Lindel\"off. 

Interesting to us, our theory in \cite{ReintjesTemple_ell2, ReintjesTemple_ell4, ReintjesTemple_ell5} couched in the language of geometry, implies a new existence theory for ordinary differential equations of type \eqref{geodesic_eqn}, based on coordinate transformation to optimal connection regularity. Building on this, we show in Section \ref{Sec_mollification} that, even though the regularity in the original $x$-coordinates is too low for the geodesic equation to admit a weak formulation (based on multiplying by smooth test functions and integrating by parts), solutions of \eqref{geodesic_eqn} in $x$-coordinates exist as limits of smooth solutions of mollified equations obtained by mollifying the connection components $\Gamma$ in a manner which maintains $L^p$ regularity of the curvature in the zero mollification limit. This shows that for the geodesic equation the notion of weak solutions in the sense of approximations is more fundamental than the notion of distributional solutions of a weak formulation of the equations, and that taking account of the geometry in the form of curvature, such weak solutions can be realized in a unique way. This raises the interesting question as to whether limits under some suitably constrained mollification might be more fundamental than weak formulations based on integration by parts, and we wonder whether this might hold some lessons for the problem of non-uniqueness for more complicated systems, like the 3-D Euler equations, c.f. \cite{DeLellisSzekelyhidi}.

Our main theorems are stated in Section \ref{Sec_results}, and their proofs are recorded in Section \ref{Sec_Proof}. In Section \ref{Sec_RT} we review authors' theory of the RT-equations, on which these proofs are based. In Section \ref{Sec_non-optimality} we characterize explicitly the loss of connection regularity relative to the curvature under coordinate transformation in terms of the exterior derivative $d\Gamma$ and its co-derivative $\delta \Gamma$.

\section{Statement of results} \label{Sec_results}

\subsection{Weak solutions by geometry}

To state our basic existence and uniqueness results, we first need to introduce a proper notion of solutions to the geodesic initial value problem \eqref{geodesic_eqn} in $x$-coordinates, given that $\Gamma$ in $L^p$ has a regularity too low to restrict to curves and hence too low for a standard weak formulation for the equations.    

\begin{Def} \label{Def_weak_geo}
We call a curve $\gamma \equiv \gamma_x$ a weak solution of \eqref{geodesic_eqn} in $x$-coordinates, if there exists a coordinate transformation $x \to y$ such that $\Gamma_x \to \Gamma_y$ transforms by the connection transformation law \eqref{connection_transfo}, $\Gamma_y \in W^{1,p},$ and the transformed curve $\gamma_y \equiv y \circ x^{-1}(\gamma_x)$ is a classical solution of \eqref{geodesic_eqn} in $y$-coordinates.  
\end{Def}

The theory of connections in differential geometry implies that weak solutions in the sense of Definition \ref{Def_weak_geo} are well-defined, because the transformation law \eqref{connection_transfo} for connections in $W^{1,p}$ preserves geometric properties under coordinate transformation, including the invariance of geodesic curves.

Our first theorem addresses existence of solutions to the geodesic initial value problem \ref{geodesic_eqn} subject to bounded curvature.

\begin{Thm} \label{Thm_existence}
Assume $\Gamma_x \in L^{2p}(\Omega)$ and ${\rm Riem}(\Gamma_x) \in L^p(\Omega)$, for $p>n$, (or equivalently, assume $\Gamma_x \in L^{2p}(\Omega)$ and $d\Gamma_x \in L^p(\Omega)$). Then there exists a  solution to the initial value problem \eqref{geodesic_eqn} in the sense of Definition \ref{Def_weak_geo}.  More generally, there exists a solution to the initial value problem for \eqref{gen_geodesic_eqn} provided $K$ is H\"older continuous.
\end{Thm}

Note that the threshold regularity for existence established by Theorem \ref{Thm_existence} lies one derivative {\it below} the H\"older continuity, $\Gamma \in W^{1,p} \subset C^{0,\alpha}$, $p>n$, $\alpha = 1 - \frac{n}{p}$, required by Peano's Theorem, if one does not take account of the curvature. Uniqueness is more delicate, and requires more regularity on ${\rm Riem}(\Gamma)$ which might differ from the regularity of $d\Gamma$, (c.f. Section \ref{Sec_non-optimality} below).

\begin{Thm} \label{Thm_uniqueness}
Assume $\Gamma_x \in L^{2p}(\Omega)$ and ${\rm Riem}(\Gamma_x) \in W^{1,p}(\Omega)$, for $p>n$. Then there exists a unique solution to the initial value problem \eqref{geodesic_eqn} in the sense of Definition \ref{Def_weak_geo}.  Moreover, if $\Gamma_x, {\rm Riem}(\Gamma_x) \in W^{1,p}(\Omega)$, $p>n$, then there exists a unique solution to the initial value problem for \eqref{gen_geodesic_eqn} provided $K$ is Lipschitz continuous.
\end{Thm}

The threshold regularity for uniqueness established by Theorem \ref{Thm_uniqueness} lies below the Lipschitz continuity, ($\Gamma \in W^{1,\infty} = C^{0,1}$), required by the Picard-Lindel\"off Theorem, if one does not take account of the curvature.  In particular, Theorem \ref{Thm_uniqueness} implies as a corollary uniqueness of solutions for H\"older continuous connections provided their curvature is in $W^{1,p}$.

To put this into the context of General Relativity (GR), note that the connection is always precisely one derivative less regular than the metric by Christoffel's formula, so the above classical thresholds for existence and uniqueness can be equivalently expressed in terms of the metric tensor, requiring $C^{1,\alpha}$ metric regularity for existence of geodesics, and $C^{1,1}$ regularity for uniqueness. Lorentzian metrics of low regularity are central to various recent research programs in GR, ranging from relativistic shock waves \cite{Israel, GroahTemple, SmollerTemple}, to Penrose's Strong Cosmic Censorship Conjecture \cite{Penrose, DafermosLuk, Kehle, Reintjes_SCC}, to the Hawking Penrose Singularity Theorems \cite{Graf, GrafGrKuSt, KunzingerOhSchSt, KunzingerStVi, HawkingEllis}, to Lorentzian Length Spaces \cite{KunzingerSaemann}. Glimm scheme based shock wave solutions in General Relativity exhibit only Lipschitz continuous gravitational metrics  \cite{Israel, GroahTemple, SmollerTemple}, i.e., one derivative below the classical threshold for existence and uniqueness of geodesics, but Theorem \ref{Thm_existence} applies since the curvature is bounded in $L^\infty$.  It was the authors' attempt to understand the apparent singularities in GR solutions generated by  Glimm's method that originally motivated the authors' to develop the theory of the RT-equations in \cite{ReintjesTemple_ell1, ReintjesTemple_ell2, ReintjesTemple_ell3, ReintjesTemple_ell4, ReintjesTemple_ell5, ReintjesTemple_ell6}, the basis for the methods in this paper.

\subsection{Weak solutions as a mollified limit}   \label{Sec_mollification}

Theorems \ref{Thm_existence} and \ref{Thm_uniqueness} use the theory of optimal regularity in geometry to define a weak solution of \eqref{geodesic_eqn} in a setting where the equation is too weak to admit either strong solutions, or solutions of a weak formulation of the equation based on integration by parts.  Having established a notion of weak solution using a higher order theory, in this case geometry, it makes sense to ask in what sense the equations are satisfied in the original $x$-coordinates, where the problem is originally posed. In this subsection we show that our notion of weak solution provides an explicit description for a mollification $\Gamma_x^{\epsilon}$ of the original connection components $\Gamma_x$, so that $\Gamma_x^{\epsilon}\to \Gamma_x$ and the weak solutions defined by Definition \ref{Def_weak_geo} are obtained as a limit of smooth solutions of the mollified equation \eqref{geodesic_eqn} in $x$-coordinates,  as $\epsilon\to 0$. (Extending the construction in this section to equation \eqref{gen_geodesic_eqn} is straightforward.)
  
To begin, note that assuming only that components of $\Gamma_x$ are functions in $L^p$, we cannot expect the weak solutions identified in Theorems \ref{Thm_existence} and \ref{Thm_uniqueness} to always be faithful limits of solutions to more regular equations obtained by arbitrarily mollifying $\Gamma_x$.   Indeed, $L^p$ functions are too weak to even restrict to curves, so in general, to realize solutions as regular limits, one would expect that mollifications need to be constrained in a manner which faithfully represents the missing physical information, in this case the transformation law for connections, and our assumptions on the curvature of $\Gamma$.   Here we prove that the RT-equations provide an explicit procedure for mollifying $\Gamma_x$ component-wise to $\Gamma_x^{\epsilon}$, so that  $\left(\Gamma_x^{\epsilon}\right)^\mu_{\rho\nu}\to (\Gamma_x)^\mu_{\rho\nu}$ in $L^{2p}$, and such that the weak solutions identified in Theorems \ref{Thm_existence} and \ref{Thm_uniqueness} are the correct limits, in $C^1$, of solutions of the regularized equation \eqref{geodesic_eqn} obtained by substituting $\Gamma_x^{\epsilon}$ for $\Gamma_x$, and solving the regularized equations in $x$-coordinates; namely, $C^1$ limits of solutions of
\beq \label{mollify1.1}
\begin{cases}
\ddot\gamma^\mu_\epsilon + \left(\Gamma_x^{\epsilon}\right)^\mu_{\rho\nu}(\gamma_\epsilon) \dot\gamma^\rho_\epsilon \dot\gamma^\nu_\epsilon =0,  \cr
\gamma_{\epsilon}^\mu(t_0)=x^\mu_0,  \cr 
\dot\gamma_{\epsilon}^\mu(t_0) = v^\mu_0.
\end{cases}
\eeq

To this end, recall that authors' existence theory for the RT-equations provides the Jacobians $J$ and transformed connection $\Gamma_y$ such that $J$ integrates to a coordinate transformation $x\to y$, and transforming $\Gamma_x$ to $y$-coordinates by the connection transformation law \eqref{connection_transfo}, yields $\Gamma_y$, a connection regular enough for the transformed geodesic equation \eqref{geodesic_eqn} to admit existence and uniqueness of solutions, according to Theorems \ref{Thm_existence} and \ref{Thm_uniqueness}. The point is now that this in turn provides a procedure for smoothly mollifying $\Gamma_x$. Namely, let $(\Gamma_y^\epsilon)^\alpha_{\beta \gamma}$ be a standard smooth mollification of $(\Gamma_y)^\alpha_{\beta \gamma}$ in $y$-coordinates, and let $x_\epsilon \equiv x_\epsilon \circ y^{-1}$ be a smooth mollification of the coordinate transformation $x \circ y^{-1}$ and let $y_\epsilon \equiv y \circ x_\epsilon^{-1}$ be its inverse, where we write $x=x_\epsilon(y)$ and $y=y_\epsilon(x)$ to distinguish the coordinates from the mollified mappings. We then introduce the mollification of $\Gamma_x$ in $x$-coordinates as
\beq \label{mollified_conn}
(\Gamma_x^\epsilon)^\mu_{\rho\nu} = \tfrac{\partial x^\mu_\epsilon}{\partial y^\alpha } \Big( \tfrac{\partial y^\beta_\epsilon}{\partial x^\rho} \tfrac{\partial y^\gamma_\epsilon}{\partial x^\nu }  \, (\Gamma_y^\epsilon)^\alpha_{\beta \gamma}   +  \tfrac{\partial }{\partial x^\rho}  \tfrac{\partial y^\alpha_\epsilon}{\partial x^\nu } \Big),
\eeq
all components in \eqref{mollified_conn} expressed as function in the original $x$-coordinates. This construction underlies the following theorem.

\begin{Thm} \label{Thm_uniqueness_refined}
Assume $\Gamma_x \in L^{2p}(\Omega)$ and ${\rm Riem}(\Gamma_x) \in W^{1,p}(\Omega)$, for $p>n$.   Then the sequence of $C^\infty$ curves $\gamma_\epsilon$ which solve \eqref{mollify1.1} converge to the correct weak solution in the following sense:   
\begin{enumerate}[(i)]
\item The $C^\infty$ curves $\gamma_\epsilon$ which solve \eqref{mollify1.1} are defined on some common open interval $I \subset \R$, for each $\epsilon > 0$, where $t_0 \in I$.
\item The sequence of connections $\Gamma_x^\epsilon$ converges  to $\Gamma_x$ strongly in $L^{2p}$ and $\Riem(\Gamma_x^\epsilon)$ converges strongly to $\Riem(\Gamma_x)$ in $L^p$ as $\epsilon \to 0$.           
\item The curves $\gamma_\epsilon$ converge to $\gamma$ strongly in $C^{1}$ as $\epsilon \to 0$, where $\gamma$ is the unique weak solution identified in Theorem \ref{Thm_uniqueness}.   
\end{enumerate}
Moreover, for \eqref{gen_geodesic_eqn}, if $\Gamma_x,\; {\rm Riem}(\Gamma_x) \in W^{1,p}(\Omega)$, $p>n$, and $K$ is Lipschitz continuous, then (i) - (iii) hold for $C^\infty$ solutions $\gamma_\epsilon$ of \eqref{mollify1.1} with a standard mollification of $K$ on the right hand side.     
\end{Thm}

To summarize, Theorems \ref{Thm_uniqueness} and \ref{Thm_existence} use the higher order theory of geometry to identify the correct weak solutions of (\ref{geodesic_eqn}) in a setting in which there does not exist a standard weak formulation of the equations; and Theorems \ref{Thm_uniqueness_refined} and \ref{Thm_existence_refined} characterize these weak solutions as limits of solutions of mollified equations.   Thus, in a physical example,   limits under mollification are more fundamental than weak formulations of the equations based on integration by parts. Authors find this interesting in light of the non-uniqueness of weak solutions of the compressible Euler equations based on integration by parts, demonstrated in \cite{DeLellisSzekelyhidi}, considering that the physically correct weak solutions of compressible Euler should be zero-viscosity limits of the regularizing Navier-Stokes equations.  That is, Theorems \ref{Thm_uniqueness_refined} and \ref{Thm_existence_refined} provide the ``physically correct'' weak solutions as mollified limits via an explicit mollification procedure,  (like Navier-Stokes to Euler), in a setting where a weak formulation of the equations does not even exist.

The next theorem shows that the above mollification procedure gives rise to sequences of curves converging in $C^1$ to weak solutions in the sense of Definition \ref{Def_weak_geo}, under the weaker assumption of Theorem \ref{Thm_existence}. 

\begin{Thm} \label{Thm_existence_refined}
Assume $\Gamma \in L^{2p}(\Omega)$ and ${\rm Riem}(\Gamma) \in L^p(\Omega)$, for $p>n$, (or equivalently, assume $\Gamma \in L^{2p}(\Omega)$ and $d\Gamma \in L^p(\Omega)$). Then (i) and (ii) of Theorem \ref{Thm_uniqueness_refined} hold, and limits of subsequences of solutions $\gamma_{\epsilon}(t)$ to \eqref{mollify1.1}, as $\epsilon \to 0$, are weak solutions of \eqref{geodesic_eqn} in the sense of Definition \ref{Def_weak_geo}.    More generally, this extends to \eqref{gen_geodesic_eqn} provided  $K$ is  H\"older continuous.
\end{Thm}     

Theorems \ref{Thm_existence_refined} and \ref{Thm_uniqueness_refined} raise a larger mathematical question, namely,  whether any mollification which appropriately takes account of the curvature bounds required for Theorems \ref{Thm_existence} and \ref{Thm_uniqueness}, would yield in the limit $\epsilon\to 0$ the correct weak solution in the sense of Definition \ref{Def_weak_geo}.

\section{Non-invariance of connection regularity relative to curvature} \label{Sec_non-optimality}

We now clarify the mechanism, identified in Section \ref{Sec_intro}, by which the regularity of connections can range by coordinate transformation from one derivative above to one derivative below the regularity of the Riemann curvature, due to the difference between the transformation laws for tensors versus connections, in terms of the exterior derivative $d$ and its co-derivative $\delta$. To start, assume a connection of optimal regularity with components $\Gamma_y \in W^{2,p}$ and curvature $\Riem(\Gamma_y) \in W^{1,p}$ in $y$-coordinates. Consider now the effect of a coordinate transformation $y \to x$ with Jacobian $J^{-1} = \frac{\partial x}{\partial y}$ having regularity identical to that of the curvature, i.e., $J^{-1} \in W^{1,p}$. (We use here $J^{-1}$ for $y \to x$ to be consistent with the notation in Sections \ref{Sec_RT} - \ref{Sec_Proof}, where we consider transformations from $x \to y$ with Jacobian $J$). For ease, assume $p>n$, so that $W^{1,p}$ is closed under multiplication by Morrey's inequality. Now write the connection transformation law \eqref{connection_transfo} as
\beq \label{connection_transfo_2}
\Gamma_x = \Gammati + J^{-1} dJ ,
\hspace{1cm} \text{where} \hspace{1cm}
\Gammati^\mu_{\rho\nu}  \equiv  (J^{-1})^\mu_\alpha  J^\beta_\rho J^\gamma_\nu  \; (\Gamma_y)^\alpha_{\beta \gamma} ,
\eeq
where $(dJ)^\alpha_{\rho \nu} \equiv \tfrac{\partial }{\partial x^\rho}  \tfrac{\partial y^\alpha}{\partial x^\nu}$, and we view $\Gamma_x$ and $\Gammati$ as matrix valued $1$-forms, e.g., $\Gammati \equiv \Gammati^\mu_{\nu j} dx^j$ with matrix indices $\mu,\nu$. Then Lemma 3.3 in \cite{ReintjesTemple_ell1} implies 
\beq \label{co-deriv_Gamma}
\delta \Gamma_x = \delta \Gammati +  \langle dJ^{-1} ; dJ \rangle + J^{-1} \Delta J ,
\eeq
where $\delta$ is the co-derivative based on the Euclidean metric in $x$-coordinates and $\Delta \equiv \delta d + d \delta$ is the standard Laplacian in $x$-coordinates, and $\langle \cdot \, ; \cdot \rangle$ is a matrix valued inner product, (see \cite{ReintjesTemple_ell1} for precise definitions). The point we would like to make, now, is that \eqref{co-deriv_Gamma} implies that the co-derivative $\delta \Gamma_x$ has the regularity of $\Delta J$, thus lies in general only in $W^{-1,p}$. However, in contrast, even though $\Riem(\Gamma)$ involves derivatives of $\Gamma$, by \eqref{Riemann_transfo} the Riemann tensor transforms by contraction with undifferentiated Jacobians $J, J^{-1} \in W^{1,p}$, (as proven in Appendix \ref{Sec_weak_curvature} for distributional curvature), and this preserves its $W^{1,p}$ regularity. This establishes that connections are mapped in general from one derivative of regularity above, to one derivative below the curvature, under coordinate transformation with Jacobians at the regularity of the curvature.

Interestingly the exterior derivative $d\Gamma$, the leading order part of $\Riem(\Gamma)$, works differently than $\delta \Gamma$ due to a cancellation of second order Jacobian derivatives. This is the reason why $d\Gamma$ can be taken in place of $\Riem(\Gamma)$ in Theorem \ref{Thm_existence}, but cannot be taken in place of $\Riem(\Gamma)$ in Theorem \ref{Thm_uniqueness}. To see this, note first that 
\beq  \label{dGamma}
d\Gamma_x = Curl(\Gamma_x) \equiv 
\tfrac{\partial}{\partial x^\tau} \Gamma^\mu_{\nu\rho} - \tfrac{\partial}{\partial x^\rho}\Gamma^\mu_{\nu\tau}
\eeq
is the leading order part of the Riemann curvature tensor   
\beq \label{Riem}
{\rm Riem}(\Gamma) = d\Gamma + \Gamma \wedge \Gamma
\eeq 
both in $x$- and $y$-coordinates, where $\Gamma \wedge \Gamma \equiv \Gamma^\mu_{\sigma \rho} \Gamma^\sigma_{\nu \tau} - \Gamma^\mu_{\sigma\tau} \Gamma^\sigma_{\nu\rho}$ is the wedge product. By Lemma 6.1 in \cite{ReintjesTemple_ell1}, we find from \eqref{connection_transfo_2} that
\beq
d\Gamma_x = d \Gammati +d J^{-1} \wedge dJ ,
\eeq
since $d^2 J =0$. That is, the exterior derivative $d\Gamma$ only contains first order derivatives of $J$ and $J^{-1}$, and thus maintains regularity when $J, J^{-1} \in W^{2,p}$, but looses one derivative when $J, J^{-1} \in W^{1,p}$.\footnote{In the latter case, Jacobian derivatives in the transformed wedge product $\Gamma \wedge \Gamma$ cancel precisely the Jacobian derivative terms in $d\Gamma$, by which the Riemann curvature transforms as a tensor and maintains its regularity.} We conclude that under singular coordinate transformation, the regularity of $\delta\Gamma$ can be one derivative below the regularity of $d\Gamma$, which in turn can be one derivative below ${\rm Riem}(\Gamma)$. But by authors theory of the RT-equations, this can always be reversed, and $\Gamma$ can always be lifted to one derivative above ${\rm Riem}(\Gamma)$ by coordinate transformation.

\section{The RT-equations and Optimal Regularity} \label{Sec_RT}

Authors proved in \cite{ReintjesTemple_ell4} that solutions of the RT-equations furnish coordinate transformations which regularize connections to one derivative of regularity above their Riemann curvature (optimal connection regularity), c.f. Theorem 2.1 in \cite{ReintjesTemple_ell4}. By this, the RT-equations extend optimal regularity and Uhlenbeck compactness to arbitrary affine connections, when before it was only known for positive definite metric geometries \cite{KazdanDeTurck, Uhlenbeck}. We now briefly review how the RT-equations establish optimal regularity, referring to \cite{ReintjesTemple_ell1, ReintjesTemple_ell4} for detailed definitions and proofs. The RT-equations are derived in \cite{ReintjesTemple_ell1} from the connection transformation law, guided by the Riemann-flat condition in \cite{ReintjesTemple_geo}. A simplified version of the original RT-equations is obtained by making a serendipitous gauge-type transformation, which uncouples the equations for the regularizing Jacobian $J\equiv \frac{\partial y}{\partial x}$ from the equations for the connection of optimal regularity, leading to what we call in \cite{ReintjesTemple_ell4} the {\it reduced} RT-equations,
\begin{eqnarray} 
\Delta J &=& \delta ( J \mm \Gamma ) - B ,  \label{RT_red_J} \\
d \vec{B} &=& \overrightarrow{\text{div}} \big(dJ \wedge \Gamma\big) + \overrightarrow{\text{div}} \big( J\, d\Gamma\big) , \label{RT_red_B1}   \\
\delta \vec{B} &=& w, \label{RT_red_B2}
\end{eqnarray}
where we view $J$,$B$ and $\Gamma \equiv \Gamma_x$ as matrix valued differential forms with components expressed in $x$-coordinates. In \cite{ReintjesTemple_ell4}, assuming $\Gamma_x \in L^{2p}$, $d\Gamma_x \in L^p$ in $x$-coordinates, we proved existence of solutions $(J,B)$ of the reduced RT-equations, $J\in W^{1,2p}$, $B\in L^{2p}$, with $J$ invertible and integrable to coordinates, by an explicit iteration scheme.  To show $J$ transforms $\Gamma_x$ to optimal regularity, we introduce the associated ``connection field'' $\Gammati$ by
\beq \label{Gammati}
\Gammati \equiv \Gamma_x - J^{-1} dJ.
\eeq
It is then proven, by exact cancellation of uncontrolled derivative terms $\delta\Gamma_x$, that $\Gammati$ solves the ``gauge transformed'' {\it first} RT-equation\footnote{The {\it reduced} RT-equations \eqref{RT_red_J} - \eqref{RT_red_B2} together with the {\it first} RT-equation \eqref{RT_first_B} are equivalent to the original {\it RT-equations}, c.f. \cite{ReintjesTemple_ell6}.} 
\beq \label{RT_first_B}
\Delta \Gammati = \delta d \Gamma_x - \delta \big( dJ^{-1} \wedge dJ \big) + d\big(J^{-1}A \big),
\eeq
where $A\equiv (B - \langle d J ; \tilde{\Gamma}\rangle)$, from which we infer $\Gammati \in W^{1,p}$ by elliptic regularity theory \cite{ReintjesTemple_ell4}. Moreover, integration of the Jacobian $J\equiv \frac{\partial y}{\partial x}$ yields a coordinate transformation $x \to y$. In light of \eqref{Gammati}, the connection  in $y$-coordinates given by
\beq \label{Gamma_y_reverse}
(\Gamma_y)^\gamma_{\alpha\beta} = J_k^\gamma (J^{-1})^i_\alpha  (J^{-1})^j_\beta \; \Gammati^k_{ij} 
\eeq
is of optimal regularity, $\Gamma_y \in W^{1,p}$. All this is proven in full detail and at the adequate level of weak solutions in \cite{ReintjesTemple_ell4},  starting from the assumption that $\Gamma \in L^{2p}$ and $d\Gamma \in L^p$ in $x$-coordinates, an assumption equivalent to $\Gamma \in L^{2p}$ and ${\rm Riem}(\Gamma) \in L^p$ by \eqref{dGamma} - \eqref{Riem}; see \cite{ReintjesTemple_ell6} for a non-technical summary.   

The generality of the setting addressed in this paper is possible because the RT-equations themselves apply in such generality--requiring nothing other than the connection components given locally in a coordinate system, making no symmetry assumptions, no requirement of a metric, nor any other technical assumptions about the background geometry, other than the regularity of the components of $\Gamma$ and $d\Gamma$ in a coordinate system--assuming no more than what is required to formulate the problem of optimal regularity.

\section{Proofs of the theorems}  \label{Sec_Proof}

\subsection{Existence - Proof of Theorem \ref{Thm_existence}}

The idea of proof is to use the RT-equations to construct a coordinate transformation $x\to y$ which regularizes the connection $\Gamma$  to H\"older continuity, the threshold regularity required for existence of solutions to \eqref{geodesic_eqn} in $y$-coordinates by the Peano Theorem. So let $\Gamma_x$ denote the coefficients of a connection in $x$-coordinates on some open and bounded set $\Omega_x \subset \R^n$. For convenience, we view $\Gamma_x$ as a connection represented in $x$-coordinates in some coordinate chart $(x,\Omega)$, on some $n$-dimensional manifold $\M$ with $\Omega_x \equiv x(\Omega) \subset \R^n$. (To reiterate, the global structure of $\M$ is not relevant here.)   Assume $\Gamma_x \in L^{2p}(\Omega_x)$ and ${\rm Riem}(\Gamma_x) \in L^p(\Omega_x)$, for $p>n$. By \eqref{dGamma} - \eqref{Riem}, this is equivalent to $\Gamma_x \in L^{2p}(\Omega_x)$ and $d\Gamma_x \in L^p(\Omega_x)$, the incoming assumption of authors'  optimal regularity result \cite[Thm 2.1]{ReintjesTemple_ell4}. Let $Q\in \Omega$ be the point where initial data is assigned in \eqref{geodesic_eqn}, i.e., $\gamma(t_0)=Q$. 

By Theorem 2.1 in \cite{ReintjesTemple_ell4}, there exists a neighborhood $\Omega' \subset \Omega$ of $Q$ on which a $W^{2,2p}$ coordinate transformation $x \to y$ is defined such that the connection components $\Gamma_y$ in $y$-coordinates have optimal regularity, $\Gamma_y \in W^{1,p}(\Omega_y')$. By Morrey's inequality \cite{Evans}, $W^{1,p}(\Omega_y') \subset C^{0,\alpha}(\Omega_y')$ for $p>n$ and $\alpha = 1 - \frac{n}{p}$, (after the usual potential change on a set of measure zero). This implies that $\Gamma_y$ is H\"older continuous. By Peano's Theorem \cite[Thm 2.1]{Hartman}, existence of at least one $C^{2}$ solution $\gamma_y$ to the initial value problem of the geodesic equation \eqref{geodesic_eqn} in $y$-coordinate now follows. Transforming the resulting geodesic curve $\gamma_y(t)$ back to $x$-coordinates with the inverse $W^{2,2p}$ coordinate transformation $y\to x$, transforming $\dot\gamma^\mu_y(t)$ as a vector and $\Gamma_y$ as a connection, the transformed curve $\gamma_x \equiv x \circ y^{-1} (\gamma_y)$ is a weak solution of \eqref{geodesic_eqn} in $x$-coordinates by Definition \ref{Def_weak_geo}. This completes the proof of existence. 

To extend the proof to the general equation \eqref{gen_geodesic_eqn}, note that tensor transformation of the vector field $K^\mu$ by $J, J^{-1} \in W^{1,2p} \subset C^{0,\alpha}$ preserves the H\"older continuity of $K^\mu$. Thus the above argument to prove existence of solutions to \eqref{geodesic_eqn} applies to \eqref{gen_geodesic_eqn} unchanged, completing the proof.
\hfill $\Box$

\subsection{Uniqueness - Proof of Theorem \ref{Thm_uniqueness}}

The idea of proof is to construct a coordinate transformation which regularizes the connection $\Gamma$  to Lipschitz continuity, the threshold regularity required by the Picard-Lindel\"off Theorem for uniqueness, by using the RT-equations twice. Let $\gamma(t_0)=Q \in \Omega$, and assume $\Gamma_x \in L^{2p}(\Omega_x)$ and $\Riem(\Gamma_x) \in W^{1,p}(\Omega_x)$, for $p>n$. By \eqref{Riem}, this implies $\Gamma_x \in L^{2p}(\Omega_x)$ and $d\Gamma_x \in L^p(\Omega_x)$, the incoming assumption of the optimal regularity result \cite[Thm 2.1]{ReintjesTemple_ell4}. By \cite[Thm 2.1]{ReintjesTemple_ell4}, there now exists a neighborhood $\Omega' \subset \Omega$ of $Q$ on which a coordinate transformation $x \to y'$ is defined, such that in $y'$-coordinates $\Gamma_{y'} \in W^{1,p}(\Omega_y')$, and the Jacobian $J'$ of the regularizing transformation has regularity $J' \in W^{1,2p}(\Omega_y')$. 

It follows by the transformation law for the curvature \eqref{Riemann_transfo} that $\Riem(\Gamma_{y'})$ maintains its $W^{1,p}$ regularity, because contraction by Jacobians in $W^{1,2p}$ does not lower its regularity, since $W^{1,p}$ is closed under multiplication by Morrey's inequality for $p>n$. Thus, in $y'$-coordinates we have $\Gamma_{y'} \in W^{1,p}(\Omega_{y'})$ and $\Riem(\Gamma_{y'}) \in W^{1,p}(\Omega_{y'})$ for $p>n$, which is the starting assumption of authors' prior optimal regularity result \cite[Thm 1.1]{ReintjesTemple_ell2}. By Theorem 1.1 in \cite{ReintjesTemple_ell2}, there exists another neighborhood $\Omega'' \subset \Omega'$ of $Q$ on which another coordinate transformation $y' \to y''$ is defined with Jacobian $J'' \in W^{2,p}$, such that in $y''$-coordinates $\Gamma_{y''} \in W^{2,p}(\Omega_y'')$.  Now, since $p>n$, Morrey's inequality implies that $W^{2,p}(\Omega_y'') \subset C^{1,\alpha}(\Omega_y'')$, and $C^{1,\alpha}(\Omega_y'') \subset W^{1,\infty} \simeq C^{0,1}$. Thus $\Gamma_{y''} \in W^{2,p}(\Omega_y'')$ is Lipschitz continuous and the Picard-Lindel\"off Theorem implies the existence of a unique $C^{2}$ solution $\gamma_{y''}$ to the initial value problem of \eqref{geodesic_eqn} in $y''$-coordinates, \cite[Thm 1.1]{Hartman}. 

Transformation of $\gamma_{y''}$ back to $x$-coordinates gives a weak solution in the sense of Definition \ref{Def_weak_geo}. Moreover, this is the only such weak solution of \eqref{geodesic_eqn} in $x$-coordinates. Namely, given a curve $\gamma_{y'''}$ which solves the transformed initial value problem \eqref{geodesic_eqn} in another coordinate system $y'''$ with $\Gamma_{y'''}$ in $W^{1,p}$, then transforming from $y'''$- to $y''$-coordinates, would regularize $\Gamma_{y'''}$ from $W^{1,p}$ to $\Gamma_{y''} \in W^{2,p}(\Omega_y'')$, and hence take $\gamma_{y'''}$ to $\gamma_{y''}$ by uniqueness of solutions to the initial value problem in $y''$-coordinates. This proves uniqueness of solutions to \eqref{geodesic_eqn}.

To prove uniqueness of solutions to \eqref{gen_geodesic_eqn}, taking into account that we assume $\Gamma$ in $W^{1,p}$, we only need to apply the second step in the above argument for proving uniqueness of solutions to \eqref{geodesic_eqn}. That is, the regularization by the coordinate transformation $y' \to y''$ suffices to prove uniqueness, because tensor transformation of the vector field $K^\mu$ by $J'', (J'')^{-1} \in W^{2,p} \subset C^{0,1}$ preserves the Lipschitz continuity of $K^\mu$.\footnote{Lipschitz continuity of $K$ might not be preserved under the $W^{1,p}$ Jacobian of the coordinate transformation from $x \to y'$, which is the reason why we need the stronger assumption $\Gamma_x \in W^{1,p}$.} This completes the proof.
\hfill $\Box$

\subsection{Existence under mollification - Proof of Theorem \ref{Thm_existence_refined}}

To begin, consider the regularized connection $\Gamma_y \in W^{1,p}(\Omega_y)$ constructed in the proof of Theorem \ref{Thm_existence}. Introduce a standard mollifier $\Gamma^\epsilon_y$ of $\Gamma_y$ which, by construction, converges strongly to $\Gamma_y$ in $W^{1,p}$, c.f. \cite{Evans}. By Morrey's inequality $\Gamma^\epsilon_y$ also converges to $\Gamma_y$ in $C^{0,\alpha}$, $\alpha = 1 - \frac{n}{p}$. Moreover, $\Riem(\Gamma^\epsilon_y)$ converges to $\Riem(\Gamma_y)$ in $L^p$. Namely, using H\"older's inequality, 
\begin{align} \label{mollified_curvature_convergence}
\|& \Riem(\Gamma_y^\epsilon)  - \Riem(\Gamma_y) \|_{L^p} 
\leq  \|d\Gamma_y^\epsilon - d\Gamma_y \|_{L^p} + \|\Gamma_y^\epsilon \wedge \Gamma_y^\epsilon - \Gamma_y \wedge \Gamma_y \|_{L^p} \cr
&\leq   \|d\Gamma_y^\epsilon - d\Gamma_y \|_{L^p} + \big( \|\Gamma_y^\epsilon\|_{L^{2p}} + \|\Gamma_y\|_{L^{2p}} \big)  \|\Gamma_y^\epsilon - \Gamma_y \|_{L^{2p}}   \cr
& \leq \big(1 + \|\Gamma_y^\epsilon\|_{L^{2p}} + \|\Gamma_y\|_{L^{2p}}  \big) \|\Gamma_y^\epsilon - \Gamma_y \|_{W^{1,p}}
\ \  \overset{\epsilon \to 0}{\longrightarrow} 0.
\end{align} 

Now there exists a unique solution $\gamma_y^\epsilon$ to the geodesic initial value problem \eqref{geodesic_eqn} for each $\Gamma_y^\epsilon$ with $\gamma_y^\epsilon \in C^{\infty}(I_\epsilon,\Omega_y)$, on open intervals $I_\epsilon \subset \R$ containing the initial time $t_0$. The convergence $\Gamma^\epsilon_y \to \Gamma_y$ in $C^{0,\alpha}$ implies $\Gamma^\epsilon_y$ is uniformly bounded  in terms of $\| \Gamma_y\|_{W^{1,p}}$ in $C^{0,\alpha}$, by the Morrey inequality, i.e.,
\beq \label{estimate_C1alpha_0}
\|\Gamma^\epsilon_y\|_{C^{0,\alpha}} 
\leq C \|\Gamma^\epsilon_y\|_{W^{1,p}} 
\leq \| \Gamma_y\|_{W^{1,p}}.
\eeq
The uniform bound \eqref{estimate_C1alpha_0} in turn implies that there exists a common open subinterval $I \subset I_\epsilon$ for all $\epsilon >0$ which contains the initial time $t_0$, and on which the solutions $\gamma_y^\epsilon$ are defined, (c.f. Appendix \ref{Sec_appendix}). Moreover, expressing the geodesic equation in \eqref{geodesic_eqn} as a first order system, standard ODE theory implies the uniform bound  
\beq \label{estimate_C1alpha}
\|\gamma_y^\epsilon\|_{C^{0,\alpha}} + \|\dot\gamma_y^\epsilon\|_{C^{0,\alpha}} \leq C
\eeq
where $C>0$ is a constant depending only on the uniform bound on $\|\Gamma^\epsilon_y\|_{C^{0,\alpha}}$ in \eqref{estimate_C1alpha_0}, the initial data and the domain $I$, c.f. \eqref{ODE_est_gamma} in the appendix. 

By \eqref{estimate_C1alpha}, the Arzela Ascoli Theorem implies $C^0$-convergence of a subsequence of $(\gamma_y^\epsilon, \dot\gamma_y^\epsilon)$. That is, a subsequence of $\gamma_y^\epsilon$ converges in $C^1$ to some limit curve $\gamma_y$. Since $\Gamma_y^\epsilon$ converges to $\Gamma_y$ in $C^{0,\alpha}$, it follows further that $\gamma_y$ is a weak solution of the initial value problem \eqref{geodesic_eqn} in $y$-coordinates in the standard sense. That is, for every test functions $\phi \in C^\infty_0(I,\R)$ the following limit holds
\beq \nonumber
\int_I \big(- \dot\gamma_y \; \dot\phi  + \phi \; \Gamma_y(\gamma_y) \dot\gamma_y  \dot\gamma_y  \big) dt   
= \lim\limits_{\epsilon \to 0} \int_I \big(- \dot\gamma^\epsilon_y \; \dot\phi  +  \phi \;\Gamma^\epsilon_y(\gamma^\epsilon_y) \dot\gamma^\epsilon_y  \dot\gamma^\epsilon_y  \big) dt  =0,
\eeq
where we omit indices and write $\Gamma_y(\gamma_y) \dot\gamma_y  \dot\gamma_y$ in place of $(\Gamma_y)^\mu_{\rho\nu}(\gamma_y) \dot\gamma_y^\rho  \dot\gamma_y^\nu$. To see this, estimate closeness of each term separately, in particular, note that 
\small $$
|\Gamma_y^\epsilon(\gamma_y^\epsilon) - \Gamma_y(\gamma_y)| 
\leq |\Gamma_y^\epsilon(\gamma_y^\epsilon) - \Gamma_y(\gamma_y^\epsilon)| + |\Gamma_y(\gamma_y^\epsilon) - \Gamma_y(\gamma_y)| 
\leq \|\Gamma_y^\epsilon - \Gamma_y\|_{C^0} + C \|\gamma_y^\epsilon - \gamma_y\|_{C^0}^\alpha
$$    \normalsize
by H\"older continuity of $\Gamma$.  Moreover, since $\Gamma_y(\gamma_y)$ and $\dot{\gamma}_y$ are both continuous, the standard weak form of the geodesic equation in $y$-coordinates implies that the weak derivative of the $C^0$ curve $\dot{\gamma}_y$ is continuous, which implies by the theory of distribution that $\gamma_y\in C^2(I)$. Thus $\gamma_y$ is in fact a classical strong solution of \eqref{geodesic_eqn} in $y$-coordinates.

We now show that the transformed curve $\gamma_x = x\circ y^{-1} (\gamma_y)$ is indeed the zero mollification limit of solutions to \eqref{mollify1.1} and thus a weak solution of \eqref{geodesic_eqn} in the sense of Definition \ref{Def_weak_geo}, as asserted by Theorem \ref{Thm_existence_refined}.  For this, we map the sequence of connections $\Gamma_y^\epsilon$ and geodesics $\gamma_y^\epsilon$ to $x$-coordinates. However, care must be taken since the coordinate transformation $y \to x$ is only in $W^{2,2p}$, and would hence not maintain smoothness of $\Gamma_y^\epsilon$ as required in Theorems \ref{Thm_uniqueness_refined} and \ref{Thm_existence_refined}. To circumvent this problem, we mollify the coordinate transformation $x\circ y^{-1}$, producing mappings $x_\epsilon \circ y^{-1}$ with Jacobians $J^{-1}_\epsilon \equiv \frac{\partial x_\epsilon}{\partial y} \in C^\infty$ and $J_\epsilon \equiv \frac{\partial y}{\partial x_\epsilon} \in C^\infty$ converging in $W^{1,2p}$ to $J^{-1}$ and $J$ respectively. Now, $x_\epsilon \circ y^{-1}$ maps each $\Gamma_y^\epsilon \in W^{1,p}$ to some $\Gamma_x^\epsilon \in C^\infty$, while maintaining $L^{2p}$ closeness of connections under the connection transformation law, 
\beq
\|\Gamma_x^\epsilon - \Gamma_x \|_{L^{2p}} \leq C^3 \big(\|\Gamma_y^\epsilon - \Gamma_y \|_{L^{2p}} + \|dJ_\epsilon - dJ \|_{L^{2p}} \big)  \overset{\epsilon \to 0}{\longrightarrow} 0  .
\eeq
Here we view each $\Gamma_x^\epsilon$ as expressed in $x$-coordinates and base norms in $x$-coordinates, we use Morrey's inequality to bound $L^\infty$ norms on Jacobians, and $C>1$ denotes a constant bounding the $W^{1,2p}$-norm of $J^{-1}_\epsilon$ and $J_\epsilon$ uniformly. Likewise, the closeness of the curvature established in \eqref{mollified_curvature_convergence} is maintained by the tensor transformation law,
\beq
\| \Riem(\Gamma_x^\epsilon)  - \Riem(\Gamma_x) \|_{L^p} 
\leq C^4 \| \Riem(\Gamma_y^\epsilon)  - \Riem(\Gamma_y) \|_{L^p}   \overset{\epsilon \to 0}{\longrightarrow} 0 .
\eeq
Moreover, $C^1$ convergence of the geodesics $\gamma_y^\epsilon$ is preserved under tensor transformation of $\dot{\gamma}_y^\epsilon$, that is,
\begin{eqnarray} 
\|\gamma_x^\epsilon - \gamma_x \|_{C^1}  &=& \|\gamma_x^\epsilon - \gamma_x \|_{C^0} + \|\dot\gamma_x^\epsilon - \dot\gamma_x \|_{C^0}  \cr
&\leq & \|\gamma_y^\epsilon - \gamma_y \|_{C^0} + C \|\dot\gamma_y^\epsilon - \dot\gamma_y \|_{C^0}   \overset{\epsilon \to 0}{\longrightarrow} 0 .
\end{eqnarray}
In summary, keeping in mind the $W^{1,2p}$ convergence of $J_\epsilon$ and $J^{-1}_\epsilon$ to $J$ and $J^{-1}$, respectively, we have proven that $\gamma_x = x\circ y^{-1} (\gamma_y)$ is a weak solution in $x$-coordinates in the sense of Definition \ref{Def_weak_geo}, and that $\gamma_x$ is the zero mollification limit satisfying (i) - (ii) of Theorem \ref{Thm_uniqueness_refined}. This completes the proof of Theorem \ref{Thm_existence_refined} for the geodesic equation.

To extend the proof to the generalized geodesic equation \eqref{gen_geodesic_eqn}, recall from the proof of Theorem \ref{Thm_existence} that H\"older continuity of $K$ is preserved under the $W^{2,2p}$ coordinate transformation from $x \to y$ and vice versa. By the same principle, closeness of a standard mollification $K^\epsilon$ of $K$ with respect to the H\"older norm is maintained under coordinate transformation. By this, it is straightforward to adapt the above analysis to $C^\infty$ solutions $\gamma_\epsilon$ of \eqref{mollify1.1} with a standard mollification of $K$ on the right hand side, and prove the assertion of Theorem \ref{Thm_existence_refined}. This completes the proof of Theorem \ref{Thm_existence_refined}.
\hfill $\Box$

\subsection{Uniqueness under mollification - Proof of Theorem \ref{Thm_uniqueness_refined}}

Consider the unique geodesic curve $\gamma_y \in C^2(I,\Omega_y)$ in $y$-coordinates constructed in the proof of Theorem \ref{Thm_uniqueness} above, where we denote here $y''$-coordinates simply by $y$. The mollification procedure in the proof of Theorem \ref{Thm_existence_refined} yields again a $C^2$ geodesic $\gamma_y$ in the zero mollification limit, which is now identical to the unique geodesic constructed in the proof of Theorem \ref{Thm_uniqueness}. Transformation back to $x$-coordinates, following again the procedure in the proof of Theorem \ref{Thm_existence_refined}, shows that $\gamma_x = x\circ y^{-1} (\gamma_y)$ is the unique weak solution in the sense of Definition \ref{Def_weak_geo}, which satisfies (i) - (iii) of Theorem \ref{Thm_uniqueness_refined}. This completes the proof for the geodesic equation.

To prove the assertion of uniqueness of Theorem \ref{Thm_existence_refined} for the generalized geodesic equation \eqref{gen_geodesic_eqn}, under the stronger assumption of $\Gamma_x \in W^{1,p}(\Omega)$, $p>n$, recall from the proof of Theorem \ref{Thm_uniqueness} that Lipschitz continuity of $K$ is preserved under the $W^{3,p}$ coordinate transformation from $y' \to y''$ and vice versa. Adapt now the analysis to $C^\infty$ solutions $\gamma_\epsilon$ of \eqref{mollify1.1} with a standard mollification of $K$ on the right hand side, the assertion of Theorem \ref{Thm_uniqueness_refined} follow again.
\hfill $\Box$

\vspace{.2cm}
Note that, despite the stronger regularity assumption ${\rm Riem}(\Gamma) \in W^{1,p}$ in Theorem \ref{Thm_uniqueness_refined} versus Theorem \ref{Thm_existence_refined}, convergence of ${\rm Riem}(\Gamma^\epsilon)$ to ${\rm Riem}(\Gamma)$ does not hold in $W^{1,p}$ in general, since closeness of the wedge product of $\Gamma$'s in \eqref{mollified_curvature_convergence} is only controlled in $L^p$ when the connection is in $L^{2p}$, as assumed in Theorem \ref{Thm_uniqueness_refined}. Starting with the stronger assumption of $\Gamma \in W^{1,p}$, convergence of the curvature in $W^{1,p}$ could also be established by adapting the above analysis in the proof of Theorem \ref{Thm_existence_refined}.

\section{Discussion of uniqueness in the singular case $\Gamma, d\Gamma$ in $L^\infty$} \label{Sec_Discussion}

Assuming $\Gamma \in L^\infty$ and $d\Gamma\in L^\infty$, we currently do not have a proof that solutions of the RT-equations furnish a regularizing coordinate transformation $x \to y$ to optimal regularity $\Gamma_y \in W^{1,\infty} \hat{=} C^{0,1}$, the threshold regularity for uniqueness of solutions to \eqref{geodesic_eqn} by the Picard Lindel\"off Theorem. This is because $p=\infty$ is a singular case of elliptic PDE theory, so the Laplacian may fail to lift solutions two derivatives above sources in $L^\infty$, and hence may fail to lift solutions two derivatives above sources in $W^{-1,\infty}$ as well.   Thus even though the RT-equations may establish the regularity $\Gamma_y \in W^{1,p}$ for any $p<\infty$, they may fail to regularize to the threshold regularity $\Gamma_y \in W^{1,\infty} \hat{=} C^{0,1}$ required for uniqueness of geodesics. Likewise, authors' investigation of the space of functions of bounded mean oscillations (BMO), (a space larger than $L^\infty$ which is contained in all $L^p$ spaces for $p<\infty$), indicates that even though the Laplacian lifts solutions of the Poisson equation two derivatives above source functions in BMO, the (non-linear) RT-equation don't appear to do so because BMO is not closed under multiplication. For comparison, in Theorems \ref{Thm_uniqueness} and \eqref{Thm_uniqueness_refined}, we assume $\Gamma\in L^{2p}$ below $L^\infty$, but this requires assuming $d\Gamma \in W^{1,p}$, a regularity above $L^\infty$. At this stage authors do not know whether there exist $L^\infty$ connections with Riemann curvature bounded in $L^\infty$, which cannot be lifted to $W^{1,\infty}$ by coordinate transformation, and for which solutions of the geodesic equation are non-unique.

\appendix

\section{Relevant ODE theory} \label{Sec_appendix}

For completeness we now review some standard estimates for systems of ODE's
\beq \label{ODE}
\dot{u} = F(u), \hspace{1cm} u(t_0)=u_0
\eeq
where $F: \Omega \to \R^m$ is assumed to be H\"older continuous on an open and bounded domain $\Omega \subset \R^n$, with bounded H\"older norm (some $0 < \alpha \leq 1$)           
\begin{eqnarray} 
\|F\|_{C^{0,\alpha}(\Omega)} &\equiv &  \|F\|_{C^0(\Omega)} + \sup_{u_1,u_2 \in \Omega}  \tfrac{| F(u_1) - F(u_2) |}{ |u_1-u_2|^\alpha},  \label{Hoelder} \\
\|F\|_{C^0(\Omega)} &\equiv & \sup_{u \in \Omega} |F(u)|.   \label{sup-norm}
\end{eqnarray}

We now derive estimate \eqref{estimate_C1alpha}, in the proof of Theorem \eqref{Thm_existence_refined}, assuming continuity of $F$ with bounded sup-norm \eqref{sup-norm}. For this, assume solutions on some bounded interval $I$ which we assume for simplicity to have length $|I|\leq 1$. We find from \eqref{ODE} and \eqref{Hoelder} (for $t>t'$) that 
\beq \label{ODE_eqn1}
|u(t) - u(t')| \; \leq \; \int^t_{t'} |\dot u| dt  
\; \leq\;  \int^t_{t'} |F(u)| dt 
\;\leq\; \|F\|_{C^0} |t-t'|  ,
\eeq
which by $|I| \leq 1$ implies $\|u\|_{C^1(I)} \equiv \|u\|_{C^0(I)} + \|\dot{u}\|_{C^0(I)} \leq \|F\|_{C^0} + |u_0|$ and
\beq \label{ODE_Hoelder-estimate}
\|u\|_{C^{0,\alpha}(I)}   
\equiv \   \|u\|_{C^0(I)} + \sup_{t_1,t_2 \in \Omega}  \tfrac{| u(t_1) - u(t_2) |}{ |t_1-t_2|^\alpha}   
\ \leq \ \|F\|_{C^0(\Omega)} + |u_0|.
\eeq 

Analogously, given a sequence of functions $F_\epsilon$ with uniform bound $\|F_\epsilon\|_{C^{0}} <C$, \eqref{ODE_Hoelder-estimate} turns into a uniform bound on $\|u_\epsilon\|_{C^{0,\alpha}(I_\epsilon)}$ for solutions $u_\epsilon$ to \eqref{ODE} with $F_\epsilon$ in place of $F$, and with fixed initial data $u_\epsilon(t_0)=u_0$, defined on intervals $I_\epsilon$. The intervals $I_\epsilon$ depend only on the uniform bound $\|F_\epsilon\|_{C^{0}(\Omega)} <C$ and the volume of $\Omega$, (c.f. \cite[Thm 2.1]{Hartman}), which implies that there exist an open subinterval $I\subset I_\epsilon$ for all $\epsilon>0$, and we assume again without loss that $|I|\leq 1$. 

To derive \eqref{estimate_C1alpha} from \eqref{ODE_Hoelder-estimate}, first use that the geodesic equation can be written in form \eqref{ODE} with $u=\big(\gamma(t),v(t)\big)$ and $F(u)= (v, \Gamma^\mu_{\sigma\rho}(\gamma) v^\sigma v^\nu)$ for $v \equiv \dot{\gamma}$. Now, assuming $\Gamma$ is defined on some set $\Omega \subset \R^n$, and restricting the domain $I$ of solutions such that $|v-\dot{\gamma}(t_0)|\leq 1$ with respect to the Euclidean norm $|\cdot|$ on $\R^n$, \eqref{ODE_Hoelder-estimate} gives 
\begin{eqnarray} \label{ODE_est_gamma}
\|\gamma\|_{C^{0,\alpha}(I)} + \|\dot{\gamma}\|_{C^{0,\alpha}(I)} 
&\leq & \|F\|_{C^0} + |u_0|,  \cr
&\leq &  b_n \|\Gamma\|_{C^0(\Omega)} + b_n^2  + |\gamma(t_0)| + |\dot{\gamma}(t_0)|,
\end{eqnarray}
where $b_n$ denotes the volume of the ball of radius $1$ in $\R^n$ (to take account of the sup-norm of $v$ over the ball $\{|v-\dot{\gamma}(t_0)|\leq 1\}$). Replacing now $\Gamma$ and $\gamma$ by $\Gamma_\epsilon$ and $\gamma_\epsilon$, \eqref{ODE_est_gamma} implies the sought after uniform $C^{1,\alpha}$ bound \eqref{estimate_C1alpha}. 

\vspace{.2cm}
We end this section with a comparison of standard H\"older versus Lipschitz estimates and their relation to uniqueness of solutions of ODE's. Combining \eqref{ODE} and \eqref{Hoelder} we obtain the standard ODE estimate on two solutions $u_1$ and $u_2$ of \eqref{ODE},          
\beq \label{ODE_eqn2}
\frac{d}{dt} |u_1-u_2|  \leq C |u_1-u_2|^\alpha,
\eeq
where $C \equiv \|F\|_{C^{0,\alpha}(\Omega)}$. For $\alpha <1$, division by $|u_1-u_2|^\alpha$ and subsequent and integration gives the basic ODE estimate
\beq \label{ODE_eqn3}
|u_1-u_2|(t)^{1-\alpha} \leq |u_1-u_2|(t_0)^{1-\alpha} + C |t-t_0|.
\eeq
This estimate is insufficient to yield uniqueness because of the growth in $t$ on the right hand side. On the other hand, in the case of Lipschitz continuity, when $\alpha =1$, the same operation as above leads to control of $\frac{d}{dt} \ln |u_1-u_2|$, and integration yields
\beq \label{ODE_eqn4}
|u_1-u_2|(t) \leq |u_1-u_2|(t_0) e^{C |t-t_0|},
\eeq
which implies uniqueness for solutions with $u_1(t_0) =u_2(t_0)$.

\section{Weak formulations of curvature and its invariance} \label{Sec_weak_curvature}

In this section we introduce the weak form of the Riemann curvature based on the Koszul formula, $\Riem(\Gamma)=d\Gamma+\Gamma\wedge\Gamma$, and we prove that, when the curvature is in $L^p$, a weak formulation based in each coordinate system is consistent with both the transformation law for the connection, and the tensor transformation of the $L^p$ components of the curvature. This establishes a notion of invariance of the weak curvature sufficient for our purposes in this paper, in particular, the control of regularity of the curvature under coordinate transformation addressed in Section \ref{Sec_non-optimality}. 

We define the weak form of the Riemann curvature tensor of a connection $\Gamma \equiv \Gamma_x$ in $L^{2p}(\Omega)$ in $x$-coordinates component-wise as a functional over the space of test functions in the sense of distributions, 
\beq \label{weak_curvature}
{\rm Riem}(\Gamma_x)[\psi]^\mu_{\nu\rho\tau} 
\equiv - \int_\omega \Big(  \Gamma^\mu_{\nu\rho} \frac{\partial \psi}{\partial x^\tau} - \Gamma^\mu_{\nu\tau} \frac{\partial\psi }{\partial x^\rho} \Big) \;dx 
+ \int_\Omega \big( \Gamma \wedge \Gamma \big)^\mu_{\nu\rho\tau}\; \psi \;dx ,
\eeq
by shifting derivatives in $Curl(\Gamma) \equiv \tfrac{\partial}{\partial x^\tau} \Gamma^\mu_{\nu\rho} - \tfrac{\partial}{\partial x^\rho}\Gamma^\mu_{\nu\tau}$ onto test functions $\psi \in C^\infty_0(\Omega)$.   We use here scalar valued test functions and the non-invariant volume element $dx$, following the setting in \cite{SmollerTemple}. We say that the weak form of the curvature \eqref{weak_curvature} is in $L^p$, if there exist functions $\mathcal{R}^\mu_{\nu\rho\tau} \in L^p(\Omega)$ such that
\beq \label{weak_curvature_Lp}
{\rm Riem}(\Gamma)[\psi]^\mu_{\nu\rho\tau} 
= \int_\Omega \mathcal{R}^\mu_{\nu\rho\tau}  \; \psi \; dx  
\ \equiv \mathcal{R}[\psi]^\mu_{\nu\rho\tau}
\eeq
for all $\psi \in C^\infty_0(\Omega)$, in which case we write $\Riem(\Gamma) = \mathcal{R} \in L^p(\Omega)$. We now prove that the weak curvature in $L^p$ transforms as a tensor under coordinate transformation following the ideas in \cite{SmollerTemple}.

In the general setting of affine connections in this paper, where no metric is assumed, there is no invariant volume element. So one cannot expect \eqref{weak_curvature} to define the exact same functional in every coordinate system. However, for our purposes it suffices to show that transforming the $L^p$ component functions $(\mathcal{R}_x)^\mu_{\nu\rho\tau}$ of the curvature in $x$-coordinates as a tensor from $x\to y$ by \eqref{Riemann_transfo} is consistent with the weak form \eqref{weak_curvature}. That is, it suffices to prove that, when transforming $\Gamma_x$ by the connection transformation law \eqref{connection_transfo} to $\Gamma_y$, then the weak form of the curvature $\Riem(\Gamma_y)$ in $y$-coordinates, defined by replacing $x$ by $y$ everywhere in \eqref{weak_curvature}, agrees with the tensor transformed $L^p$ functions $(\mathcal{R}_y)^\delta_{\alpha \beta \gamma}$ in $y$-coordinates in the sense of \eqref{weak_curvature_Lp}. This is accomplished in the following lemma.

\begin{Lemma} \label{Lemma_weak_curvature}
Consider a coordinate transformation $x \to y$ with Jacobian $J \equiv \tfrac{\partial y}{\partial x} \in W^{1,p}(\Omega_x)$, $p>n$. Let $\Gamma_x \in L^{2p}(\Omega_x)$ and $\Gamma_y  \in L^{2p}(\Omega_y)$ be the connection components in $x$- and $y$-coordinates, respectively, related by the connection transformation law \eqref{connection_transfo}. Assume that $\Riem(\Gamma_x) = \mathcal{R}_x$ for a collection of functions $\mathcal{R}_x \in L^p(\Omega)$  in $x$-coordinates  in the sense of \eqref{weak_curvature_Lp}. Define $\mathcal{R}_y$ by the tensor transformation law \eqref{Riemann_transfo}, i.e., $(\mathcal{R}_x)^\tau_{\mu\nu \rho}=  \tfrac{\partial x^\tau }{\partial y^\delta} \; \tfrac{\partial y^\alpha}{\partial x^\mu } \tfrac{\partial y^\beta}{\partial x^\nu } \tfrac{\partial y^\gamma}{\partial x^\rho } \; (\mathcal{R}_y)^\delta_{\alpha \beta \gamma}$. Then $\mathcal{R}_y \in L^p(\Omega_y)$ and $\Riem(\Gamma_y) = \mathcal{R}_y$ in the sense of \eqref{weak_curvature_Lp}.
\end{Lemma}

\Proof
Assume for the moment the coordinate transformation $x\to y$ is smooth. By assumption, $\mathcal{R}_x$ and $\mathcal{R}_y$ are related by the tensor transformation law \eqref{Riemann_transfo},  
\beq \label{weak_1}
(\mathcal{R}_x)^\tau_{\mu\nu \rho}=  \tfrac{\partial x^\tau }{\partial y^\delta} \; \tfrac{\partial y^\alpha}{\partial x^\mu } \tfrac{\partial y^\beta}{\partial x^\nu } \tfrac{\partial y^\gamma}{\partial x^\rho } \; (\mathcal{R}_y)^\delta_{\alpha \beta \gamma},
\eeq
which we write for brevity as $\mathcal{R}_x =  \tfrac{\partial y}{\partial x} \; \mathcal{R}_y$. For our incoming assumption ${\rm Riem}(\Gamma_x) = \mathcal{R}_x \in L^p(\Omega_x)$, we obtain from \eqref{weak_curvature_Lp} under coordinate transformation
\begin{eqnarray} \label{weak_2}
{\rm Riem}(\Gamma_x)[\psi] 
= \int_{\Omega_x} \mathcal{R}_x  \; \psi \; dx 
= \int_{\Omega_y} \tfrac{\partial y}{\partial x} \; \mathcal{R}_y \; \psi \; \big| \tfrac{\partial x}{\partial y} \big| \; dy,
\end{eqnarray}
where $\big| \tfrac{\partial x}{\partial y} \big|$ denotes the determinant of the inverse Jacobian $\tfrac{\partial x}{\partial y}$ resulting by transformation of the volume element from $x\to y$, and we transform $\psi$ as a scalar function, $\psi(y) = \psi(x(y))$. 

On the other hand, for a standard mollification $\Gamma_x^\epsilon$ of the connection components $\Gamma_x$, since \eqref{weak_curvature} only involves the undifferentiated connection components in $L^{2p}(\Omega_x) \subset L^1(\Omega_x)$, we obtain
\begin{eqnarray} \label{weak_3'}
{\rm Riem}(\Gamma_x)[\psi] 
 =  \lim_{\epsilon\to 0} \; {\rm Riem}(\Gamma_x^\epsilon)[\psi]  
= \lim_{\epsilon\to 0} \int_{\Omega_x} \mathcal{R}_x^\epsilon  \; \psi \; dx  .
\end{eqnarray}
Here $\mathcal{R}_x^\epsilon$ denotes the smooth components of the curvature ${\rm Riem}(\Gamma_x^\epsilon)$, equivalent to the weak form \eqref{weak_curvature} for the mollified connections $\Gamma_x^\epsilon$. Moreover, transforming $\Gamma_x^\epsilon$ to $\Gamma_x^\epsilon$ by the connection transformation law \eqref{connection_transfo} under the smooth transformation $x\to y$, the components $\mathcal{R}_y^\epsilon$ of $ {\rm Riem}(\Gamma_y^\epsilon)$ transform by the tensor transformation law \eqref{Riemann_transfo}, which gives  
\begin{eqnarray} \label{weak_3}
{\rm Riem}(\Gamma_x)[\psi] 
= \lim_{\epsilon\to 0} \int_{\Omega_y} \tfrac{\partial y}{\partial x} \; \mathcal{R}_y^\epsilon \; \psi \; \big| \tfrac{\partial x}{\partial y} \big| \; dy  .
\end{eqnarray}

Now, since \eqref{weak_2} and \eqref{weak_3} hold for all $\psi \in C^\infty_0(\Omega)$, it follows that
\begin{eqnarray} \label{weak_4}
\int_{\Omega_y}  \mathcal{R}_y \; \psi \; dy  
 =  \lim_{\epsilon\to 0} \int_{\Omega_y}  \mathcal{R}_y^\epsilon \; \psi \; dy  .
\end{eqnarray}
By smoothness, integration by parts on the right hand side of  \eqref{weak_4} gives the weak form \eqref{weak_curvature}, where the $L^{2p}$ convergence of $\Gamma_y^\epsilon$ converges to $\Gamma_y$, (as a result of the $L^{2p}$ convergence of $\Gamma_x^\epsilon$ and the connection transformation law), yields convergence of the weak form of the curvature,
\begin{eqnarray} \label{weak_5}
\int_{\Omega_y}  \mathcal{R}_y \; \psi \; dy  
 =   \lim_{\epsilon\to 0}  \Riem(\Gamma_y^\epsilon)[\psi] 
= \Riem(\Gamma_y)[\psi].
\end{eqnarray}
This proves that $\Riem(\Gamma_y) = \mathcal{R}_y$ in the sense of \eqref{weak_curvature_Lp} under smooth coordinate transformations, while $\mathcal{R}_y \in L^p(\Omega_y)$ follows directly from the tensor transformation \eqref{weak_1}. The proof directly extends to coordinate transformations with Jacobians in $W^{1,p}$ by mollification, using that $W^{1,p}$ is closed under multiplication for $p>n$. This completes the proof. 
\QED

\section*{Funding}
M. Reintjes was partially supported by CityU Start-up Grant for New Faculty (7200748) and by CityU Strategic Research Grant (7005839).

\section*{Acknowledgement}
We thanks Craig Evans for helpful comments on elliptic regularity theory in the singular case of $L^\infty$ based Sobolev spaces.

\end{document}